\documentstyle[11pt,twoside,leqno]{article}
\input epsf
\input amssym

\def\be{\begin{eqnarray}}
\def\ee{\end{eqnarray}}

\begin{document}
\date{}
\title{Symplectic Connections Induced by the Chern Connection}
\author{ Ebrahim Esrafilian  and  Hamid Reza Salimi Moghaddam}

\maketitle \abstract{} \setcounter{equation}{0} Let $(M,\omega)$
be a symplectic manifold and $F$ be a Finsler structure on $M$. In
the present paper we define a lift of the symplectic two-form
$\omega$ on the manifold $TM\backslash 0$, and find the conditions
that the Chern connection of the Finsler structure $F$ preserves
this lift of $\omega$. In this situation if $M$ admits a nowhere
zero vector field then we have a non-empty family of Fedosov
structures on $M$.



\section{Introduction.}
\setcounter{equation}{0} Connections are important objects in
differential geometry. In 1917 a connection introduced on a
manifold embedded into $R^n$ by Levi-Civita. One year later H.
Weyl introduced general symmetric linear connections in the
tangent bundle. In 1922 E. Cartan studied non-symmetric linear
connection which was applied in general relativity as a drastic
tool(\cite{[GeReSh]}).\\
An interesting field about the connections in differential
geometry is the relation between connections and the other
structures on a manifold. One of these structures is the
symplectic structure on the manifolds. A symplectic connection is
a symmetric connection which preserves the symplectic form. Many
mathematicians worked on the symplectic connections. They found
many interesting results about the relations of symplectic
connections and the other structures on the manifolds and
supermaniflds (see \cite{[BaCa1]}, \cite{[BaCa2]},
\cite{[BiCaGuRaSch]}, \cite{[GeReSh]}, \cite{[LaRa]} and
\cite{[Sch]}.).\\
Let $(M,\omega)$ be a symplectic manifold. A natural question is:
"Which connections preserve $\omega$?".\\
Conversely, let $\nabla$ be a connection on a manifold $M$, we
can ask ourselves: "Which nondegenerete closed two-forms $\omega$
exist such that $\nabla$ preserves them?".\\
I. Gelfand, V. Retakh and M. Shubin have obtained some interesting
results about the first question in \cite{[GeReSh]}. In this
paper we try to answer to the second question in a special case.\\
Let $(M,F)$ be a Finsler manifold also at the same time let
$(M,\omega)$ be a symplectic manifold. Suppose that
$\tilde{\nabla}$ is the Chern connection arising from the Finsler
structure $F$. We define a canonical lift of $\omega$ on the
manifold $TM\backslash 0$ to find the conditions for $\omega$,
such that $\tilde{\nabla}$ preserves the canonical lift of
$\omega$. In this case if $M$ admits a nowhere zero vector field
then we have a non-empty family of symplectic connections on $M$.

\section{Preliminaries and notations.}
\setcounter{equation}{0}

In this section we give some important definitions and theorems
which we need in the following. In the first, we review some of
general definitions and theorems of symplectic geometry.

\textbf{Definition 2.1.} A \textit{symplectic form} (or a
\textit{symplectic structure}) on a manifold $M$ is a
nondegenerate, closed two-form $\omega$ on $M$. A
\textit{symplectic manifold} $(M,\omega)$ is a manifold $M$
together with a symplectic form $\omega$ on $M$ (For more details
see \cite{[AbMa]}).

\textbf{Definition 2.2.} Assume $(M,\omega)$ is a symplectic
manifold. Let $\nabla$ be a connection (covariant derivative) on
$M$. We say that $\nabla$ \textit{preserves} $\omega$ if
$\nabla\omega=0$ or \be
Z(\omega(X,Y))=\omega(\nabla_ZX,Y)+\omega(X,\nabla_ZY), \nonumber
\ee for any vector fields $X, Y, Z$ (For more details see
\cite{[GeReSh]}.).

\textbf{Definition 2.3.} If $\nabla$ is symmetric and preserve
the given symplectic form $\omega$ then we say that $\nabla$ is a
\textit{symplectic connection}.

\textbf{Definition 2.4.} \textit{Fedosov manifold} is a
symplectic manifold with a given symplectic connection. We show a
Fedosov manifold by $(M,\omega,\nabla)$, where $\nabla$ is the
symplectic connection such that preserves $\omega$.

\textbf{Theorem 2.5.} (Darboux) Suppose $\omega$ is a
nondegenerate two-form on a $2n-$\\ manifold $M$. Then $d\omega=0$
if and only if there is a chart $(U,\phi)$ at each $m\in M$ such
that $\phi(m)=0$, and with $\phi(u)=(x^1(u),\cdots,x^{2n}(u))$ we
have \be\omega|U=\sum_{i=1}^{n}dx^i\wedge dx^{n+i}.\ee (For more
details see \cite{[AbMa]}.)

Now we give some fundamental concepts of Finsler geometry.

\textbf{Definition 2.6.} A Minkowski norm on $R^n$ is a
nonnegative function ${\mathcal{F}}:R^n\rightarrow
[0,\infty)$ which has the following properties:\\

(i) ${\mathcal{F}}$ is ${\mathcal{C}}^\infty$ on $R^n \setminus \{0\}$.\\

(ii) ${{\mathcal{F}}}(\lambda y)=\lambda {{\mathcal{F}}}(y)$ for
all $\lambda > 0$
and $y\in R^n$.\\

(iii) The $n\times n$ matrix $(g_{ij})$, where $g_{ij}(y):=
[\frac{1}{2}{{\mathcal{F}}}^2]_{y^iy^j}(y)$, is
positive-\hspace*{1cm}definite at all $y\neq 0$.

\textbf{Definition 2.7.} Let $M$ be an $n-$dimensional smooth
manifold and $TM$ the tangent bundle of $M$. A function
$F:TM\rightarrow
[0,\infty)$ is called a Finsler metric if it has the following properties:\\

(i) $F$ is ${\mathcal{C}}^\infty$ on the slit tangent bundle $TM\setminus 0$.\\

(ii) For each $x\in M$, ${\mathcal{F}}_x:=F|_{T_xM}$ is a
Minkowski norm on $T_xM$.

If the Minkowski norm satisfies ${\mathcal{F}}(-y)=
{\mathcal{F}}(y)$, then one has the absolutely homogeneity
$F(\lambda y)= |\lambda| F(y)$, for any $\lambda\in R$. Every
absolutely homogeneous Minkowski norm is a norm in the sense of
functional
analysis.\\
Every Riemannian manifold $(M,g)$ by defining \be
F(x,y):=\sqrt{g_x(y,y)} \hspace*{2cm}x\in M, y\in T_xM\nonumber\ee
is a Finsler manifold.\\ Also we use of the following
relations:\be
g_{ij}&:=&(\frac{1}{2}F^2)_{y^iy^j}=FF_{y^iy^j}+F_{y^i}F_{y^j}, \nonumber\\
A_{ijk}&:=&\frac{F}{2}\frac{\partial g_{ij}}{\partial
y^k}=\frac{F}{4}(F^2)_{y^iy^jy^k}, \nonumber\\
\frac{\delta y^i}{F}&:=&\frac{1}{F}(dy^i+N^i_jdx^j), \nonumber\ee
where the $N^i_j$ are collectively known as the nonlinear
connection.

\textbf{Theorem 2.8.}  (Chern) Let $(M,F)$ be a Finsler manifold.
The pulled-back bundle $\pi^\ast TM$ admits a unique linear
connection, called the Chern connection. Its connection forms
$(\theta^i_j)$ are characterized by the structural
equations:\\

$\ast$ Torsion freeness:
$d(dx^i)-dx^j\wedge\theta^i_j=-dx^j\wedge\theta^i_j=0$,\\

$\ast$ Almost $g-$compatibility:
$dg_{ij}-g_{kj}\theta^k_i-g_{ik}\theta^k_j=2A_{ijs}\frac{\delta
y^s}{F}$.\\ (For more details see \cite{[BaChSh]}.)

Also we use of the following notations.\\
Let $(x^i)$ and $(\hat{x}^i),i=1,\cdots, n$ be local coordinates,
we use of the following notations:\\
(1) $\partial_i:=\frac{\partial}{\partial x^i}$ and
$\hat{\partial}_i:=\frac{\partial}{\partial \hat{x}^i}$.\\
(2) $\omega_{ij}:=\omega(\partial_i,\partial_j)$.\\
(3) $\nabla_i:=\nabla_{\partial_i}$.\\
(4) (\textit{Christoffel symbols} $\Gamma_{ij}^k$),
$\nabla_i\partial_j=\Gamma_{ij}^k\partial_k$ for symplectic
connections.\\
(5) (\textit{Christoffel symbols} $\tilde{\Gamma}_{ij}^k$),
$\tilde{\nabla}_i\tilde{\partial}_j=\tilde{\Gamma}_{ij}^k\tilde{\partial}_k$
for the Chern connection $\tilde{\nabla}$.

\section{Symplectic Connections Induced by the Chern Connection.}

\setcounter{equation}{0} For Christoffel symbols of a symplectic
connection we have the following theorem.

\textbf{Theorem 3.1.} \textit{Let $(M,\omega,\nabla)$ be a Fedosov
manifold. In Dorboux coordinate $(x^i)_{i=1}^{2n}$ for any
$k=1,\cdots,2n$ we have:
\be\Gamma_{kj}^{i+n}&=&\Gamma_{ki}^{j+n}\hspace*{1cm} 1\leq i\leq
n \ ,\ 1\leq j\leq n, \nonumber\\
\Gamma_{kj}^{i+n}&=&-\Gamma_{ki}^{j-n}\hspace*{1cm} 1\leq i\leq n
\ ,\ n+1\leq j\leq 2n, \nonumber\\
\Gamma_{kj}^{i-n}&=&-\Gamma_{ki}^{j+n}\hspace*{1cm} n+1\leq i\leq
2n \ ,\ 1\leq j\leq n, \nonumber\\
\Gamma_{kj}^{i-n}&=&\Gamma_{ki}^{j-n}\hspace*{1cm} n+1\leq i\leq
2n \ ,\ n+1\leq j\leq 2n. \nonumber \ee}\\

\textit{Proof}: It is suffices to use the equation
$\omega=\sum_{i=1}^n dx^i\wedge dx^{n+i}$ and use of it in the
following equation: \be
0=\partial_k\omega(\partial_i,\partial_j)=\omega(\nabla_k\partial_i,\partial_j)
+\omega(\partial_i,\nabla_k\partial_j).\nonumber \ee\hfill\ $\Box$\\

Now we define a canonical lift of forms on the manifold
$TM\backslash 0$.

\textbf{Definition 3.2.} Suppose that $\omega$ is a $k-$form on a
manifold $M$ of dimension $n$. In this case we define a $k-$form
$\tilde{\omega}$ on the manifold $TM\backslash 0$, associated to
pulled-back bundle $\pi^\ast TM$, in the following way and call
it the \textit{canonical lift of $\omega$ on $\pi^\ast TM$}: \be
\tilde{\omega}_{(x,y)}: \overbrace{(\pi^\ast
TM)_{(x,y)}\times\cdots\times (\pi^\ast
TM)_{(x,y)}}^{k-times}\rightarrow R \nonumber\\
\tilde{\omega}_{(x,y)}(z_1,\cdots,z_n)=\omega_x(z_1,\cdots,z_n).
\nonumber \ee For any $x\in M, y\in T_xM\backslash 0$ and
$z_1,\cdots,z_n\in T_xM$.

\textbf{Definition 3.3.} Let $X\in {\mathcal{X}}(M)$ be a vector
field on $M$. We define a section $\tilde{X}$ of $\pi^\ast TM$ in
the following way and call it the \textit{canonical lift of $X$}:
\be \tilde{X}_{(x,y)}=X_x. \nonumber\ee For any $x\in M$ and $y\in
T_xM\backslash 0$. We use of the notation $\partial_i$ with sense
$\tilde{\partial}_i$.

\textbf{Definition 3.4.} Let $(M,\omega)$ be a symplectic
manifold with a Finsler structure $F$ on it. In this case we call
the triple $(M,\omega,F)$ a \textit{Finslerian symplectic
manifold}.

\textbf{Definition 3.5.} Suppose that $(M,\omega,F)$ is a
Finslerian symplectic manifold and $\tilde\nabla$ is the Chern
connection of $(M,F)$. We say that $\tilde\nabla$ preserves
$\omega$ if $\tilde\nabla$ preserves $\tilde\omega$, the canonical
lift of $\omega$, in the other words if for any $X, Y, Z \in
{\mathcal{X}}(M)$ we have \be
\tilde{Z}(\tilde{\omega}(\tilde{X},\tilde{Y}))=
\tilde{\omega}(\tilde{\nabla}_{\tilde{Z}}\tilde{X},\tilde{Y})+\tilde{\omega}(\tilde{X},
\tilde{\nabla}_{\tilde{Z}}\tilde{Y}).\nonumber\ee Where
$\tilde{X}, \tilde{Y}, \tilde{Z}$ are the canonical lift of $X, Y,
Z$ and
$\tilde{Z}(\tilde{\omega}(\tilde{X},\tilde{Y})):=Z(\omega(X,Y))$.
Suppose that
$\tilde{\omega}_{ij}:=\tilde{\omega}(\tilde{\partial}_i,\tilde{\partial}_j)$
then for $X=\partial_i,Y=\partial_j,Z=\partial_k$ we have \be
\tilde{\partial}_k\tilde{\omega}_{ij}=
\tilde{\omega}_{il}\tilde{\Gamma}_{kj}^l-\tilde{\omega}_{jl}\tilde{\Gamma}_{ki}^l.\nonumber\ee

\textbf{Theorem 3.6.} Let $(M,\omega,F)$ be a Finslerian
symplectic manifold and $\tilde{\nabla}$ be the Chern connection
of $(M,F)$ such that $\tilde{\nabla}$ preserves $\omega$. Suppose
that $M$ admits a vector field $W\in{\mathcal{X}}(M)$ such that
for any $x\in M$ we have $W_x\neq 0_{T_xM}$. Then there is a
nonempty family of symplectic connections on $M$ such that
preserves $\omega$, in the other words there is a nonempty family
of Fedosov structures on $M$.

\textit{Proof}: By using Christoffel symbols
$\tilde{\Gamma}_{ij}^k$ of the Chern connection $\tilde{\nabla}$
on $TM\backslash 0$ and vector field $W$ we define new Christoffel
symbols $\Gamma_{ij}^k$ on $M$ in the following way: \be
\Gamma_{ij}^k(x):=\tilde{\Gamma}_{ij}^k(x,W_x). \nonumber\ee Now
we have a symmetric linear connection $\nabla$, generated by
Christoffel symbols $\Gamma_{ij}^k$, on $M$ such that preserves
$\omega$ because for $\partial_i, \partial_j$ and $\partial_k$ we
have: \be \{
\omega({\nabla}_{\partial_k}\partial_i,\partial_j)+{\omega}(\partial_i,
{\nabla}_{\partial_k}\partial_j)\}_x
&=&\{\omega(\Gamma_{ki}^l\partial_l,\partial_j)
+\omega(\partial_i,\Gamma_{kj}^l\partial_l)\}_x\nonumber\\
&=&\{\tilde{\omega}(\tilde{\Gamma}_{ki}^l\tilde{\partial}_l,\tilde{\partial}_j)
+\tilde{\omega}(\tilde{\partial}_i,\tilde{\Gamma}_{kj}^l\tilde{\partial}_l)\}_{(x,W_x)}\nonumber\\
&=&\{\tilde{\partial}_k\tilde{\omega}(\tilde{\partial_i},\tilde{\partial}_j)\}_{(x,W_x)}\nonumber\\
&=&\{\partial_k \omega (\partial_i,\partial_j) \}_x. \nonumber\ee
So $\nabla$ is a symplectic connection on $M$. \hfill\
$\Box$\\

\textbf{Corollary 3.7.} In Theorem 3.6, if $(M,F)$ is of Berwald
type then the induced symplectic connection is unique.

\textit{Proof}: By using the definition of Berwald spaces,
Christoffel symbols of the Chern connection are functions only of
variable $x\in M$ so the induced symplectic connection is not
associated to vector field $W$. Therefore the induced symplectic
connection is unique. \hfill\ $\Box$\\

\textbf{Theorem 3.8.} Let $(M,\omega,F)$ be a Finslerian
symplectic manifold. Suppose that $F$ is locally Minkowskian. Then
the Chern connection preserves $\omega$ if and only if \be
\partial_k\omega_{ij}&=&0,\nonumber\\
\partial_h\hat{x}^l.\hat{\partial}_k\hat{\partial}_ix^h.\omega_{lj}
&+&\partial_h\hat{x}_l.\hat{\partial}_k\hat{\partial}_jx^h.\omega_{il}
=\hat{\partial}_k\omega_{ij}. \nonumber \ee Where
$(x^i)_{i=1}^{2n}$ is the natural coordinates for the locally
Minkowskian manifold $M$ and $(\hat{x}^i)_{i=1}^{2n}$ is an
arbitrary local coordinate on $M$.

\textit{Proof}: In any Finsler manifold $M$ we have \be
\hat{\tilde{\Gamma}}_{qr}^p=\partial_i\hat{x}^p.\hat{\partial}_q\hat{\partial}_rx^i
+\partial_i\hat{x}^p.\tilde{\Gamma}_{jk}^i.\hat{\partial}_qx^j\hat{\partial}_rx^k.
\ee (See \cite{[BaChSh]} p.43) In locally Minkowskian manifolds in
natural coordinates we have $\tilde{\Gamma}_{jk}^i=0$, so by using
3.5 and equation 3.1 we obtain the result.\hfill\ $\Box$\\

Now let us use the Randers metrics to obtain the Fedosov
structure.\\
Let $F(x,y)=\alpha(x,y)+\beta(x,y)$ be a Randers metric on a
manifold $M$ such that $d\beta$ be nondegenerate. Therefore
$(M,d\beta)$ is a symplectic manifold. Now we calculate the
conditions that the Chern connection preserves $d\beta$.\\
Suppose that $\beta(x,y)=b_i(x)y^i$ or in other words
$\beta(x)=b_i(x)dx^i$ and let $\omega=d\beta$, so we have
\be\omega=d(b_idx^i)=db_i\wedge dx^i.\nonumber\ee Therefore
\be\omega_{pq}&=&(db_i\wedge
dx^i)(\partial_p,\partial_q)=\partial_pb_q-\partial_qb_p.\nonumber\ee
The Chern connection arising from $F=\alpha+\beta$, preserves
$\omega$ if and only if
\be(\tilde{\Gamma}_{ki}^l\partial_lb_j-\tilde{\Gamma}_{kj}^l\partial_lb_i)
+(\tilde{\Gamma}_{kj}^l\partial_ib_l-\tilde{\Gamma}_{ki}^l\partial_jb_l)
+(\partial_k\partial_jb_i-\partial_k\partial_ib_j)=0\nonumber\ee
and also in Darboux coordinate we have
\be(\tilde{\Gamma}_{ki}^l\partial_lb_j-\tilde{\Gamma}_{kj}^l\partial_lb_i)
+(\tilde{\Gamma}_{kj}^l\partial_ib_l-\tilde{\Gamma}_{ki}^l\partial_jb_l)=0.
\nonumber \ee

\section{Results from the Curvature Tensor}
Let $(M,\omega,F)$ be a Finslerian symplectic manifold such that
the Chern connection preserves $\omega$ and $M$ admits a vector
field $W\in{\mathcal{X}}(M)$ such that for any $x\in M$ we have
$W_x\neq 0_{T_xM}$. Suppose that $\nabla$ is the induced
symplectic connection on $M$ by the Chern connection and vector
field $W$. We have the following relations for the curvature
tensor of $\nabla$: \be
R(X,Y)Z=\nabla_X\nabla_YZ-\nabla_Y\nabla_XZ-\nabla_{[X,Y]}Z,
\nonumber \ee \be
R(\partial_j,\partial_k)\partial_i=R_{ijk}^l\partial_l. \nonumber
\ee Therefore by using the Christoffel symbols we have: \be
R_{ijk}^l(x)&=&\partial_j\Gamma_{ki}^l(x)-\partial_k\Gamma_{ij}^l(x)+\Gamma_{ki}^m(x).\Gamma_{jm}^l(x)
-\Gamma_{ij}^m(x).\Gamma_{km}^l(x)\nonumber\\
&=&\partial_j\tilde{\Gamma}_{ki}^l(x,W(x))-\partial_k\tilde{\Gamma}_{ij}^l(x,W(x))\nonumber\\
&+&\tilde{\Gamma}_{ki}^m(x,W(x)).\tilde{\Gamma}_{jm}^l(x,W(x))
-\tilde{\Gamma}_{ij}^m(x,W(x)).\tilde{\Gamma}_{km}^l(x,W(x))\nonumber\\
&=&(D_j\tilde{\Gamma}_{ki}^l)(x,W(x))+(D_p\tilde{\Gamma}_{ki}^l)(x,W(x)).\partial_jW^p(x)\nonumber\\
&-&(D_k\tilde{\Gamma}_{ij}^l)(x,W(x))-(D_p\tilde{\Gamma}_{ij}^l)(x,W(x)).\partial_kW^p(x)\nonumber\\
&+&\tilde{\Gamma}_{ki}^m(x,W(x)).\tilde{\Gamma}_{jm}^l(x,W(x))\nonumber\\
&-&\tilde{\Gamma}_{ij}^m(x,W(x)).\tilde{\Gamma}_{km}^l(x,W(x)).\label{eq1}
\ee Where in local coordinate $(x^i)$ we have \be W=W^p\partial_p,
\nonumber \ee and also $D_p$ is derivative associated to component
$p$.

\textbf{Theorem 4.1.} For any symplectic connection we have: \be
R_{ijkl}=R_{jikl}. \label{eq2} \ee\\

\textit{Proof}: See \cite{[GeReSh]}. \hfill\ $\Box$\\

\textbf{Theorem 4.2.} Let $(M,\omega,F)$ be a Finslerian
symplectic manifold. If the Chern connection $\tilde{\nabla}$
preserves $\omega$ then for any nowhere zero $W\in
{\mathcal{X}}(M)$ we have
the following two conditions:\\
(1) \be
&\{&\{D_k\tilde{\Gamma}_{lj}^n+D_p\tilde{\Gamma}_{lj}^n.\partial_kW^p
-D_l\tilde{\Gamma}_{jk}^n-D_p\tilde{\Gamma}_{jk}^n.\partial_lW^p
+\tilde{\Gamma}_{lj}^p.\tilde{\Gamma}_{kp}^n
-\tilde{\Gamma}_{jk}^p.\tilde{\Gamma}_{lp}^n\}\nonumber\\
&+&\{D_l\tilde{\Gamma}_{jk}^n+D_p\tilde{\Gamma}_{jk}^n.\partial_lW^p
-D_j\tilde{\Gamma}_{kl}^n-D_p\tilde{\Gamma}_{kl}^n.\partial_jW^p
+\tilde{\Gamma}_{jk}^p.\tilde{\Gamma}_{lp}^n
-\tilde{\Gamma}_{kl}^p.\tilde{\Gamma}_{jp}^n\}\nonumber\\
&+&\{D_j\tilde{\Gamma}_{kl}^n+D_p\tilde{\Gamma}_{kl}^n.\partial_jW^p
-D_k\tilde{\Gamma}_{lj}^n-D_p\tilde{\Gamma}_{lj}^n.\partial_kW^p
+\tilde{\Gamma}_{kl}^p.\tilde{\Gamma}_{jp}^n
-\tilde{\Gamma}_{lj}^p.\tilde{\Gamma}_{kp}^n\}\}\omega_{in}=0,\nonumber
\ee (2) \be
&\{&D_k\tilde{\Gamma}_{lj}^n+D_p\tilde{\Gamma}_{lj}^n.\partial_kW^p
-D_l\tilde{\Gamma}_{jk}^n-D_p\tilde{\Gamma}_{jk}^n.\partial_lW^p
+\tilde{\Gamma}_{lj}^p.\tilde{\Gamma}_{kp}^n
-\tilde{\Gamma}_{jk}^p.\tilde{\Gamma}_{lp}^n\}\omega_{in}\nonumber\\
&-&\{D_k\tilde{\Gamma}_{li}^n+D_p\tilde{\Gamma}_{li}^n.\partial_kW^p
-D_l\tilde{\Gamma}_{ik}^n-D_p\tilde{\Gamma}_{ik}^n.\partial_lW^p
+\tilde{\Gamma}_{li}^p.\tilde{\Gamma}_{kp}^n
-\tilde{\Gamma}_{ik}^p.\tilde{\Gamma}_{lp}^n
\}\omega_{jn}=0.\nonumber \ee\\

\textit{Proof}: For (1) it is suffices to use equation \ref{eq1}
and the first Bianchi identity for symmetric connections. We can
obtain (2) by using equations \ref{eq1} and \ref{eq2}.\hfill\ $\Box$\\

\vspace{5mm}
\begin{flushleft}
\medskip Ebrahim Esrafilian, Faculty of mathematics,
Department of pure mathematics, Iran university of science and
technology, Narmak, Tehran 16846-13114, Iran. \\

\vspace{5mm}
\medskip Hamid Reza Salimi Moghaddam, Department of Mathematics, Faculty of  Sciences, University of Isfahan, Isfahan,81746-73441-Iran. \\

\vspace{2mm} e-mail:\\
 salimi.moghaddam@gmail.com and hr.salimi@sci.ui.ac.ir
\end{flushleft}


\begin{thebibliography}{9}

\bibitem{[AbMa]} R. Abraham, J. E. Marsden, {\em Foundations of
Mechanics,} Second Edition, Addison-Wesley Publication Company,
(1987).
\bibitem{[BaCa1]}P. Baguis, M. Cahen, {A construction of symplectic connections through reduction,
\em Lett. Math. Phys., }{\bf 57} (2001), 149 - 160.
\bibitem{[BaCa2]}{P. Baguis, M. Cahen, Marsden-Weinstein reduction for symplectic connections,
\em Bull. Belg. Math. Soc. Simon Stevin 10,}{\bf No.1} (2003), 91
- 100.
\bibitem{[BaChSh]}D. Bao, S. S. Chern, Z. Shen, {\em An Introduction to Riemann-Finsler
Geometry,} Springer-Verlag, (2000).
\bibitem{[BiCaGuRaSch]}P. Bieliavsky, M. Cahen, S. Gutt, J. Rawnsley and L. Schwachh\"ofer,
{\em Symplectic connections}, Int. J. Geom. Methods Mod. Phys., 3,
no. 3, (2006), 375 - 420.
\bibitem{[GeReSh]}I. Gelfand, V. Retakh, M. Shubin, {\em Fedosov Manifolds},
Symplectic Geometry Workshop, Toronto, (June 1997).
\bibitem{[LaRa]} P. M. Lavrov, O. V. Radchenko, {\em On higher order relations in Fedosov supermanifolds},
 J. Phys. A: Math. Gen., 39 (2006), 6501 - 6508.
\bibitem{[Sch]} L. J. Schwachh\"ofer, {\em Special connections on symplectic manifolds},
Proceedings of the 24th Winter School ``Geometry and Physics'',
January 2004, Srni, Czech Republic, in: Supplemento ai Rendiconti
del Circolo Matematico di Palermo, Serie II, Nr. 75 (2005), 197 -
223.

\end{thebibliography}
\end{document}